# HYBRID DEAD-BEAT OBSERVERS FOR A CLASS OF NONLINEAR SYSTEMS

**Iasson Karafyllis[*] and Zhong-Ping Jiang[**]**

[*]Dept. of Environmental Eng., Technical University of Crete,
73100, Chania, Greece, email: ikarafyl@enveng.tuc.gr

[**]Dept. of Electrical and Computer Eng., Polytechnic Institute of New York University,
Six Metrotech Center, Brooklyn, NY 11201, U.S.A., email: zjiang@control.poly.edu

**Abstract**
This paper studies the reduced-order or full-order, dead-beat observer problem for a class of nonlinear systems, linear in the unmeasured states. A novel hybrid observer design strategy is proposed, with the help of the notion of strong observability in finite time. The proposed methodology is applied to a batch reactor, for which a hybrid dead-beat observer is obtained in the absence of the precise measurements of the concentration variables. Moreover, the observer is used for the estimation of the frequency of a sinusoidal signal. The results show that accurate estimations can be provided even if the signal is corrupted by high frequency noise.

## 1. Introduction

The observer problem occupies an important place in mathematical control theory. It is concerned with the estimation of unmeasured states of a dynamic control system using the information of inputs and outputs. There is a vast literature on the problems of existence and design of observers (see for instance [1,2,4,5,6,7,10,11,13,16,17,20,27,28,31,32] and references therein). In this work we focus on nonlinear systems of the form:

$$\dot{x}(t) = A(y(t), u(t))x(t) + b(y(t), u(t))$$
$$\dot{y}_i(t) = f_i(y(t), u(t)) + \sum_{j=1}^{n} c_{i,j}(y(t))x_j(t) \quad , \quad i = 1, ..., k \qquad (1.1)$$
$$x(t) = (x_1(t), ..., x_n(t))' \in \Re^n , \; y(t) = (y_1(t), ..., y_k(t))' \in \Re^k$$
$$(x(t), y(t)) \in O \subseteq \Re^n \times \Re^k , \; u(t) \in U \subseteq \Re^m$$

where $O \subseteq \Re^{n+k}$ is an open set, $U \subseteq \Re^m$ is a non-empty closed set, $A(y,u) = \{a_{i,j}(y,u), i, j = 1, ..., n\}$ and all mappings $a_{i,j} : \Omega \times U \to \Re$ ($i, j = 1, ..., n$), $b : \Omega \times U \to \Re^n$, $c_{i,j} : \Omega \to \Re$ ($i = 1, ..., k$, $j = 1, ..., n$) and $f_i : \Omega \times U \to \Re$ ($i = 1, ..., k$) are locally Lipschitz, where $\Omega = \{y \in \Re^k : \exists x \text{ such that } (x, y) \in O\}$. It is assumed that the component of the state vector $y$, also known as the output, is available for the feedback design and that the remaining state component $x$ is unmeasured and is to be estimated.

Systems of the form (1.1) are termed as "systems linear in the unmeasured state components" in the literature (see [3,8,9,22,23,24]). The dynamic output feedback stabilization problem has been studied extensively in the past for this class of systems in [3,8,9,22,23,24]. Exponential observers for systems linear in the unmeasured state components were provided in [4], under a persistency of excitation condition. It should be noted that systems of the form (1.1) are related to systems with output dependent incremental rate. For systems with output dependent incremental rate the dynamic high-gain approach was exploited in [25] for the solution of the output feedback stabilization problem.

The purpose of the present work is to study the observability properties of systems linear in the unmeasured state components and to propose a novel observer design procedure that guarantees features which cannot be provided by conventional observers: we propose hybrid observers which provide exact estimation of the unmeasured state components in finite time (dead-beat observers). Moreover, we consider the general case where the system evolves in an open set $O$ and not in $\Re^{n+k}$: this generality allows us to study interesting applications (see Section 3). It should be noted that hybrid observers were recently proposed in [5] as well. Moreover, dead-beat observers have been proposed in the literature for linear systems:

- by means of sliding modes (see [11,17,18]),

- by means of delays (see [7]).

The approach of using delays for the observer design was exploited in [20] for a special class of nonlinear systems (with nonlinear output injection terms). High-gain techniques were utilized in [6,28] for the design of semi-global finite-time observers for a class of nonlinear systems.

In Section 2 of the present work, we study the observability properties of system (1.1). The notion of strong observability is introduced for general nonlinear systems. Conditions for the construction of inputs which do not distinguish between events in finite time (see [30]) for system (1.1) are proposed. Then, it is shown that under strong observability in time $r > 0$ for system (1.1), it is possible to construct a deterministic system ($\Sigma$) with states $z(t) \in O$ and inputs $y : \Re^+ \to \Omega$, $u : \Re^+ \to U$ so that

$$z(t) = (x(t), y(t)), \text{ for all } t \geq r \qquad (1.2)$$

The proposed observer ($\Sigma$) is a hybrid system which uses delays: the history of the output is utilized in order to estimate the state component $x$ of system (1.1). For the overall system (1.1) with ($\Sigma$) the classical semigroup property does not hold: however, the weak semigroup property holds (see [14,15]). Also, the overall system (1.1) with ($\Sigma$) is autonomous in the sense described in [14,15]. The proposed hybrid observer relies on the minimization of an appropriate $L^2$ norm and is methodologically close to the procedure used for optimization-based observers (see [1] and references therein).

Finally, in Section 3 of the present work applications are presented for the proposed hybrid, dead-beat observer:

1. The first application was inspired by the example presented in [27]. The application deals with the operation of a batch reactor, where the temperature is continuously measured. However, in practice it is often very difficult to obtain accurate measurements of the concentrations of the reactants. Moreover, accurate estimates of the concentrations of the reactants are needed in short time, because such systems are not left to operate for large times (see [19]). The application of the proposed hybrid, dead-beat observer can provide exact estimates of the concentrations of the reactants in finite time.

2. The second application deals with the estimation of the frequency of a sinusoidal signal. The problem was recently studied in [26] (see also [12,21]). It is shown that the proposed hybrid, dead-beat observer provides robust estimation of the frequency of the measured signal. Indeed, the observer is tested in the presence of high frequency noise and exactly the same test of robustness with the one in [26] is performed. The results show that the sensitivity to measurement noise decreases as the time horizon of the minimized $L^2$ norm increases, i.e., as the length of the history of the output which is utilized for the state estimation increases. This feature is expected and it is common to optimization-based observers (see [1] and references therein).

**Notations** Throughout this paper we adopt the following notations:

∗ Let $I \subseteq \Re_+ := [0, +\infty)$ be an interval. By $L^\infty(I;U)$ ($L^\infty_{loc}(I;U)$) we denote the space of measurable and (locally) essentially bounded functions $u(\cdot)$ defined on $I$ and taking values in $U \subseteq \Re^m$.

∗ By $C^0(A;\Omega)$, we denote the class of continuous functions on $A$, which take values in $\Omega$.

∗ For a vector $x \in \Re^n$ we denote by $x'$ its transpose and by $|x|$ its Euclidean norm. The determinant of a square matrix $A \in \Re^{n \times n}$ is denoted by $\det(A)$. $A' \in \Re^{n \times m}$ denotes the transpose of the matrix $A \in \Re^{m \times n}$.

∗ By $A = diag(l_1, l_2, ..., l_n)$ we mean a diagonal matrix with $l_1, l_2, ..., l_n$ on its diagonal.

## 2. Hybrid Dead-Beat Observer Design

Consider an autonomous system described by ordinary differential equations of the form:

$$\dot{x}(t) = F(x(t), u(t))$$
$$x(t) \in D \subseteq \Re^n, \, u(t) \in U \subseteq \Re^m \tag{2.1}$$

where $D \subseteq \Re^n$ is an open set, $U \subseteq \Re^m$ is a non-empty closed set and the mapping $F : D \times U \to \Re^n$ is locally Lipschitz. The output of system (2.1) is given by

$$y(t) = h(x(t)) \tag{2.2}$$

where the mapping $h : D \to \Re^k$ is continuous. For system (2.1) we adopt the following notion of observability. We assume that for every $x_0 \in D$ and $u \in L^\infty_{loc}(\Re_+; U)$ there exists a unique solution $[0, +\infty) \ni t \to x(t) = x(t, x_0; u) \in D$ for (2.1) satisfying (2.1) for almost every $t \geq 0$ and $x(0) = x(0, x_0; u) = x_0$.

**Definition 2.1:** *Consider system (2.1) with output (2.2). We say that the input $u \in L^\infty([0, r]; U)$ strongly distinguishes the state $x_0 \in D$ in time $r > 0$, if the following condition holds*

$$\max_{t \in [0,r]} |h(x(t, x_0; u)) - h(x(t, \xi; u))| > 0, \text{ for all } \xi \in D \text{ with } x_0 \neq \xi \tag{2.3}$$

**Remark 2.2:** It is clear that Definition 2.1 implies that if the input $u \in L^\infty([0, r]; U)$ strongly distinguishes the state $x_0 \in D$ in time $r > 0$ then for every $\xi \in D$ with $x_0 \neq \xi$ the input $u \in L^\infty([0, r]; U)$ distinguishes between the events $(x_0, 0)$ and $(\xi, 0)$ (see [30]).

For system (1.1) we assume that for every $(x_0, y_0) \in O$ and $u \in L^\infty_{loc}(\Re_+; U)$ there exists a unique solution $[0, +\infty) \ni t \to (x(t), y(t)) = (x(t, x_0, y_0; u), y(t, x_0, y_0; u)) \in O$ for (1.1) satisfying (1.1) for almost every $t \geq 0$ and $(x(0), y(0)) = (x(0, x_0, y_0; u), y(0, x_0, y_0; u)) = (x_0, y_0)$.

Denote by $\Phi(t, x_0, y_0; u)$ the transition matrix of the linear time-varying system $\dot{x}(t) = A(y(t), u(t))x(t)$ corresponding to given inputs $u \in L^\infty_{loc}(\Re_+; U)$ and $y(t) = y(t, x_0, y_0; u)$ for $t \geq 0$. Then the following fact holds for the solutions of system (1.1). It follows directly from integration of the differential equations (1.1).

<u>Fact I:</u> For every $(x_0, y_0) \in O$ and $u \in L^\infty([0, r]; U)$ the following equations hold for all $t \geq 0$:

$$x(t, x_0, y_0; u) = \Phi(t, x_0, y_0; u)x_0 + \theta(t, x_0, y_0; u) \tag{2.4}$$

$$p(t, x_0, y_0; u) = q'(t, x_0, y_0; u) \, x_0 \tag{2.5}$$

where

$$q(t, x_0, y_0; u) := \int_0^t \Phi'(s, x_0, y_0; u) C(s, x_0, y_0; u) ds \tag{2.6}$$

$$\theta(t, x_0, y_0; u) := \int_0^t \Phi(t, x_0, y_0; u) \Phi^{-1}(\tau, x_0, y_0; u) b(y(\tau, x_0, y_0; u), u(\tau)) d\tau \tag{2.7}$$

$$C'(t, x_0, y_0; u) := \begin{bmatrix} c_{1,1}(y(t, x_0, y_0; u)) & \cdots & c_{1,n}(y(t, x_0, y_0; u)) \\ \vdots & & \vdots \\ c_{k,1}(y(t, x_0, y_0; u)) & \cdots & c_{k,n}(y(t, x_0, y_0; u)) \end{bmatrix} \in \Re^{k \times n} \tag{2.8}$$

$$p(t,x_0,y_0;u) := y(t,x_0,y_0;u) - y_0 - \begin{bmatrix} \int_0^t f_1(y(s,x_0,y_0;u),u(s))ds \\ \vdots \\ \int_0^t f_k(y(s,x_0,y_0;u),u(s))ds \end{bmatrix} - \int_0^t C'(s,x_0,y_0;u)\theta(s,x_0,y_0;u)ds \quad (2.9)$$

It is important to note at this point that all expressions involved in (2.4)-(2.9) can be evaluated by means of the output trajectory $y(\tau) = y(\tau, x_0, y_0; u)$ for $\tau \in [0,t]$ and the input $u(\tau)$ for $\tau \in [0,t]$. For example, the transition matrix $\Phi(t, x_0, y_0; u)$ can be evaluated by solving the linear matrix differential equation $\frac{d}{d\tau}\Phi(\tau) = A(y(\tau),u(\tau))\Phi(\tau)$ for $\tau \in [0,t]$ with initial condition $\Phi(0) = I$, where $I$ denotes the identity matrix. Similarly, $C(\tau) := C(\tau, x_0, y_0; u)$ is simply $C'(\tau) := \{c_{i,j}(y(\tau)), i=1,...,k, j=1,...,n\}$ for $\tau \in [0,t]$ and $\theta(t) = \theta(t, x_0, y_0; u)$ can be computed by solving the linear system of differential equations $\frac{d}{d\tau}\theta(\tau) = A(y(\tau),u(\tau))\theta(\tau) + b(y(\tau),u(\tau))$ for $\tau \in [0,t]$ with initial condition $\theta(0) = 0 \in \Re^n$. Finally, the differential equations $\frac{d}{d\tau}q(\tau) = \Phi'(\tau)C(\tau)$ and $\frac{d}{d\tau}\xi(\tau) = (f_1(y(\tau),u(\tau)),...,f_k(y(\tau),u(\tau)))' + C'(\tau)\theta(\tau)$, for $\tau \in [0,t]$ can be utilized to provide the quantities $q(t) = q(t, x_0, y_0; u)$ and $p(t, x_0, y_0; u) = y(t) - y(0) - \xi(t)$.

The following proposition provides characterizations of the class of inputs $u \in L^\infty([0,r];U)$ which strongly distinguish the state $(x_0, y_0) \in O$ in time $r > 0$ for system (1.1). The basic idea of Proposition 2.3 is the conversion of the observability property to the minimization of an appropriate $L^2$ norm. Therefore our approach is close to the procedure used for optimization-based observers (see [1] and references therein). The proof of the following proposition is postponed to the Appendix.

**Proposition 2.3:** *Consider system (1.1). The following statements are equivalent:*

**(a)** *The input $u \in L^\infty([0,r];U)$ strongly distinguishes the state $(x_0, y_0) \in O$ in time $r > 0$.*

**(b)** *The problem*

$$\min_{\xi \in B(y_0)} \int_0^r |p(t,x_0,y_0;u) - q'(t,x_0,y_0;u)\xi|^2 dt \quad (2.10)$$

*where $B(y_0) := \{\xi \in \Re^n : (\xi, y_0) \in O\}$, admits the unique solution $\xi = x_0$.*

**(c)** *The symmetric matrix*

$$Q(r,x_0,y_0;u) := \int_0^r q(t,x_0,y_0;u)q'(t,x_0,y_0;u)dt \quad (2.11)$$

*is positive definite. Moreover, it holds that*

$$x_0 = Q^{-1}(r,x_0,y_0;u)\int_0^r q(t,x_0,y_0;u)p(t,x_0,y_0;u)dt \quad (2.12)$$

**(d)** *The following implication holds:*

$$q'(t,x_0,y_0;u)\xi = 0 \quad, \forall t \in [0,r] \Rightarrow \xi = 0 \in \Re^n \quad (2.13)$$

**Remark 2.4:** Suppose that $k=1$ and that the input $u \in L^{\infty}([0,r];U)$ does not strongly distinguish the state $(x_0, y_0) \in O$ in time $r > 0$ for system (1.1). The equivalence (a) $\Leftrightarrow$ (d) shows that there exists $\xi \in \Re^n$, $\xi \neq 0$ such that $q'(t, x_0, y_0; u) \xi = 0$ for all $t \in [0, r]$. It follows that $\frac{d}{dt} q'(t, x_0, y_0; u) \xi = C'(t, x_0, y_0; u) \Phi(t, x_0, y_0; u) \xi = 0$, for all $t \in [0, r]$. Consequently, we obtain:

$$\det\left(\begin{bmatrix} C'(t, x_0, y_0; u) \Phi(t, x_0, y_0; u) \\ C'(t_1, x_0, y_0; u) \Phi(t_1, x_0, y_0; u) \\ \vdots \\ C'(t_{n-1}, x_0, y_0; u) \Phi(t_{n-1}, x_0, y_0; u) \end{bmatrix}\right) = 0 \text{, for all } t, t_1, \ldots, t_{n-1} \in [0, r] \quad (2.14)$$

Thus, we obtain the following corollary.

**Corollary 2.5:** *Consider system (1.1) with $k=1$ and let $(x_0, y_0) \in O$, $u \in L^{\infty}([0,r];U)$ for which there exist $t, t_1, \ldots, t_{n-1} \in [0, r]$ such that:*

$$\det\left(\begin{bmatrix} C'(t, x_0, y_0; u) \Phi(t, x_0, y_0; u) \\ C'(t_1, x_0, y_0; u) \Phi(t_1, x_0, y_0; u) \\ \vdots \\ C'(t_{n-1}, x_0, y_0; u) \Phi(t_{n-1}, x_0, y_0; u) \end{bmatrix}\right) \neq 0 \quad (2.15)$$

*Then the input $u \in L^{\infty}([0,r];U)$ strongly distinguishes the state $(x_0, y_0) \in O$ in time $r > 0$. Moreover, the symmetric matrix $Q(r, x_0, y_0; u)$ defined by (2.11) is positive definite and (2.12) holds.*

It is convenient to exploit condition (2.14) in order to inputs which do not strongly distinguish the state $(x_0, y_0) \in O$ in time $r > 0$. The following example illustrates the use of (2.14).

**Example 2.6:** Consider the system

$$\begin{aligned}
\dot{x}_1(t) &= a_1(y(t))x_1(t) \\
\dot{x}_2(t) &= a_2(y(t))x_2(t) \\
\dot{y}(t) &= u(t) + c_1(y(t))x_1(t) + c_2(y(t))x_2(t) \\
x(t) &= (x_1(t), x_2(t)) \in \Re^2, \ y(t) \in \Re, \ u(t) \in \Re
\end{aligned} \quad (2.16)$$

where $c_i : \Re \to (0, +\infty)$, $a_i : \Re \to \Re$, $i = 1, 2$ are continuously differentiable functions. For system (2.16) we assume forward completeness for all $u \in L^{\infty}_{loc}(\Re_+; \Re)$: for example, forward completeness can be guaranteed if the functions $a_i : \Re \to \Re$, $i = 1, 2$ are bounded from above and the functions $c_i : \Re \to (0, +\infty)$, $i = 1, 2$ satisfy a linear growth condition, i.e., $|c_i(y)| \leq A|y| + B$, $i = 1, 2$ for all $y \in \Re$ and for certain constants $A, B \geq 0$. By virtue of Remark 2.4, if the input $u \in L^{\infty}([0,r]; \Re)$ does not strongly distinguish the state $(x_0, y_0) \in \Re^3$ in time $r > 0$ then the following condition must hold:

$$c_2(y_0) c_1(y(t, x_0, y_0; u)) \exp\left(\int_0^t (a_1(y(s, x_0, y_0; u)) - a_2(y(s, x_0, y_0; u))) ds\right) = c_1(y_0) c_2(y(t, x_0, y_0; u)),$$
$$\text{for all } t \in [0, r] \quad (2.17)$$

The reader should notice that in this case we have $\Phi(t, x_0, y_0; u) = diag\left(\exp\left(\int_0^t a_1(y(s, x_0, y_0; u)) ds\right), \exp\left(\int_0^t a_2(y(s, x_0, y_0; u)) ds\right)\right)$. Condition (2.17) coincides with condition (2.14) for $n = 2$, $t_1 = 0$. By differentiating (2.17) we obtain:

$$\kappa(y(t,x_0,y_0;u))\dot{y}(t,x_0,y_0;u) + (a_1(y(t,x_0,y_0;u)) - a_2(y(t,x_0,y_0;u))) = 0, \text{ for almost all } t \in [0,r] \quad (2.18)$$

where

$$\kappa(y) = \frac{d}{dy}\ln\left(\frac{c_1(y)}{c_2(y)}\right) \quad (2.19)$$

Moreover, using (2.16) and (2.17) with $x_0 = (x_{1,0}, x_{2,0})' \in \Re^2$, we obtain:

$$\dot{y}(t,x_0,y_0;u) = u(t) + c_2(y(t,x_0,y_0;u))\exp\left(\int_0^t a_2(y(s,x_0,y_0;u))ds\right)\left(x_{2,0} + x_{1,0}\frac{c_1(y_0)}{c_2(y_0)}\right),$$

for almost all $t \in [0,r]$ \quad (2.20)

If we further assume that $\kappa(y) \neq 0$, for all $y \in \Re$, then we conclude from (2.18), (2.20) that:

*"If the input $u \in L^\infty([0,r];\Re)$ does not strongly distinguish the state $(x_0, y_0) \in \Re^2 \times \Re$ in time $r > 0$ with $x_0 = (x_{1,0}, x_{2,0})' \in \Re^2$ for system (2.16), then the input $u \in L^\infty([0,r];\Re)$ satisfies for almost all $t \in [0,r]$:*

$$u(t) = \frac{a_2(y(t)) - a_1(y(t))}{\kappa(y(t))} - c_2(y(t))\exp\left(\int_0^t a_2(y(s))ds\right)\left(x_{2,0} + x_{1,0}\frac{c_1(y_0)}{c_2(y_0)}\right) \quad (2.21)$$

*where $y:[0,r] \to \Re$ is the solution of the initial value problem*

$$\dot{y}(t) = \frac{a_2(y(t)) - a_1(y(t))}{\kappa(y(t))}, \text{ with } y(0) = y_0 \text{"} \quad (2.22)$$

Therefore, condition (2.14) allowed us to construct inputs which do not strongly distinguish the state $(x_0, y_0) \in \Re^3$ in time $r > 0$. Indeed, without additional hypotheses we cannot be sure that every input $u \in L^\infty([0,r];\Re)$ will strongly distinguish every state $(x_0, y_0) \in \Re^3$ in time $r > 0$. For example, if there exists $y^* \in \Re$ such that $a_1(y^*) = a_2(y^*) = 0$ then the input $u(t) \equiv u^* = -c_2(y^*)x_{2,0} - c_1(y^*)x_{1,0}$ cannot distinguish between the state $(x_0, y^*) = (x_{1,0}, x_{2,0}, y^*) \in \Re^2 \times \Re$ and the state $(\xi, y^*) = (\xi_1, \xi_2, y^*) \in \Re^2 \times \Re$ with $\xi_2 = x_{2,0} + \frac{c_1(y^*)}{c_2(y^*)}(x_{1,0} - \xi_1)$: both states produce the same output response $y(t) \equiv y^*$ when the constant input $u(t) \equiv u^* = -c_2(y^*)x_{2,0} - c_1(y^*)x_{1,0}$ is applied. ◁

We next define the notion of strongly observable systems in time $r > 0$.

**Definition 2.7:** *Consider system (2.1) with output (2.2). We say that (2.1) is strongly observable in time $r > 0$ if every input $u \in L^\infty([0,r];U)$ strongly distinguishes every state $x_0 \in D$ in time $r > 0$.*

**Remark 2.8:** Proposition 2.3 guarantees that system (1.1) is strongly observable in time $r > 0$, if the symmetric matrix $Q(r,x_0,y_0;u)$ defined by (2.11) is positive definite for all $(x_0,y_0) \in O$ and $u \in L^\infty([0,r];U)$. It is clear that observability for the linear time-invariant system $\dot{x} = Ax + Bu$, $y = c'x$, where $x \in \Re^n$, $y \in \Re^k$, $A \in \Re^{n \times n}$, $B \in \Re^{n \times m}$, $c \in \Re^{k \times n}$ is equivalent to strong observability in time $r > 0$ for every $r > 0$. In general, strong observability in time $r > 0$ implies observability in time $r > 0$ in the sense of [30]. However, for nonlinear systems the converse statement does not hold. Notice that system (2.16) of Example 2.6 is not strongly observable in time $r > 0$; however, it is observable in time $r > 0$ in the sense described in [30]: every input $u \in L^\infty([0,r];\Re)$ which does not satisfy (2.21) for almost all $t \in [0,r]$ is an input which distinguishes between the events $(x_0, y_0, 0)$ and $(\xi, y_0, 0)$ in time $r > 0$, where $(x_0, y_0) = (x_{1,0}, x_{2,0}, y_0) \in \Re^3$, $(\xi, y_1) = (\xi_1, \xi_2, y_1) \in \Re^3$.

Proposition 2.3 shows that, under the following hypotheses for system (1.1):

**(H1)** *System (1.1) is strongly observable in time $r > 0$.*

then we are in a position to define the operator:

$$P : C^0([0,r];\Omega) \times L^\infty([0,r];U) \to \Re^n$$

where $\Omega = \{y \in \Re^k : \exists x \text{ such that } (x, y) \in O\}$. For each $(y, u) \in C^0([0,r];\Omega) \times L^\infty([0,r];U)$, $P(y, u)$ is defined by

$$P(y,u) = \Phi(r, y; u) Q^{-1} \int_0^r q(\tau) p(\tau) d\tau + \theta(r) \tag{2.23}$$

where $\Phi(t, y; u)$ is the transition matrix of the linear system $\dot{z}(t) = A(y(t), u(t)) z(t)$, $Q = \int_0^r q(\tau) q'(\tau) d\tau$,

$q(\tau) = \int_0^\tau \Phi'(s, y; u) C(s) ds$, $C'(\tau) := \{c_{i,j}(y(\tau)), i = 1,...,k, j = 1,...,n\}$, $p(\tau) = y(\tau) - y(0) - \int_0^\tau f(y(s), u(s)) ds - \int_0^\tau C'(s) \theta(s) ds$,

$f(y, u) := (f_1(y, u), ..., f_k(y, u))'$, $\theta(\tau) := \int_0^\tau \Phi(\tau, y; u) \Phi^{-1}(s, y; u) b(y(s), u(s)) ds$ for all $\tau \in [0, r]$. Proposition 2.3 guarantees that, if hypothesis (H1) holds for system (1.1), then for every $(x_0, y_0) \in O$ and $u \in L^\infty_{loc}(\Re_+; U)$ the following equality holds:

$$x(t, x_0, y_0; u) = P(\delta_{t-r} y, \delta_{t-r} u), \text{ for all } t \geq r \tag{2.24}$$

where $(\delta_{t-r} y)(s) = y(t - r + s, x_0, y_0; u)$, $(\delta_{t-r} u)(s) = u(t - r + s)$ for $s \in [0, r]$.

Therefore, if hypothesis (H1) holds for system (1.1), then we are in a position to provide a hybrid, dead-beat observer for system (1.1). Given $t_0 \geq 0$, $(z_0, w_0) \in O$, we calculate $(z(t), w(t)) \in O$ by the following algorithm:

Step $i$ : Calculation of $z(t)$ for $t \in [t_0 + ir, t_0 + (i+1)r]$

1) Calculate $z(t)$ for $t \in [t_0 + ir, t_0 + (i+1)r]$ as the solution of $\dot{z}(t) = A(w(t), u(t)) z(t) + b(w(t), u(t))$, $\dot{w}_i(t) = f_i(w(t), u(t)) + \sum_{j=1}^n c_{i,j}(w(t)) z_j(t)$ ($i = 1, ..., k$), with $w(t) = (w_1(t), ..., w_k(t))' \in \Re^k$.

2) Set $z(t_0 + (i+1)r) = P(\delta_{t_0 + ir} y, \delta_{t_0 + ir} u)$ and $w(t_0 + (i+1)r) = y(t_0 + (i+1)r)$, where $P : C^0([0,r];\Omega) \times L^\infty([0,r];U) \to \Re^n$ is the operator defined by (2.23).

For $i = 0$ we take $(z(t_0), w(t_0)) = (z_0, w_0)$ (initial condition).

The proposed observer can be represented by the following system of equations:

$$\dot{z}(t) = A(w(t), u(t))z(t) + b(w(t), u(t)) \, , \, t \in [\tau_i, \tau_{i+1})$$

$$\dot{w}_i(t) = f_i(w(t), u(t)) + \sum_{j=1}^{n} c_{i,j}(w(t))x_j(t) \quad , \quad i = 1,\ldots,k \, , \, t \in [\tau_i, \tau_{i+1})$$

$$z(\tau_{i+1}) = P(\delta_{\tau_i} y, \delta_{\tau_i} u)$$

$$w(\tau_{i+1}) = y(\tau_{i+1}) \tag{2.25}$$

$$\tau_{i+1} = \tau_i + r$$

$$z(t) = (z_1(t),\ldots,z_n(t))' \in \Re^n \, , \, w(t) = (w_1(t),\ldots,w_k(t))' \in \Re^k$$

$$(z(t), w(t)) \in O \subseteq \Re^n \times \Re^k$$

Thus, from all the above results, we obtain the following corollary.

**Corollary 2.9:** *Consider system (1.1) and assume that hypothesis (H1) holds. Consider the unique solution $(x(t), y(t), z(t), w(t)) \in O \times O$ of (1.1), (2.25) with arbitrary initial condition $(x_0, y_0, z_0, w_0) \in O \times O$ corresponding to arbitrary input $u \in L^\infty_{loc}(\Re_+; U)$. Then the solution $(x(t), y(t), z(t), w(t)) \in O \times O$ of (1.1), (2.25) satisfies:*

$$z(t) = x(t) \text{ and } w(t) = y(t), \text{ for all } t \geq r \tag{2.26}$$

**Remark 2.10:** The proposed observer (2.25) is a hybrid system which uses delays: the history of the output is utilized in order to estimate the state component $x$ of system (1.1). For the overall system (1.1) with (2.25) the classical semigroup property does not hold: however, the weak semigroup property holds (see [14,15]). Also, the overall system (1.1) with (2.25) is autonomous in the sense described in [14,15]. Finally, it should be noted that by virtue of Corollary 2.5 for the case $k = 1$ a sufficient condition for hypothesis (H1) is the following condition:

*"For every $(x_0, y_0) \in D \times \Omega$, $u \in L^\infty([0,r];U)$ there exists $t, t_1,\ldots,t_{n-1} \in [0,r]$ such that (2.15) holds"*

Next assume that the following hypothesis holds in addition to hypothesis (H1).

**(H2)** *There exist open sets $D \subseteq \Re^n$ and $\Omega \subseteq \Re^k$ such that $O = D \times \Omega$. Moreover, for every $\xi \in D$ and for every $(y, u) \in C^0([0,r];\Omega) \times L^\infty([0,r];U)$, the solution $z(t)$ of $\dot{z}(t) = A(y(t), u(t))z(t) + b(y(t), u(t))$ with $z(0) = \xi$ satisfies $z(t) \in D$ for all $t \in [0,r]$.*

If hypothesis (H2) holds then we can design a reduced-order, hybrid, dead-beat observer for system (1.1) of the form:

$$\dot{z}(t) = A(y(t), u(t))z(t) + b(y(t), u(t)) \, , \, t \in [\tau_i, \tau_{i+1})$$

$$z(\tau_{i+1}) = P(\delta_{\tau_i} y, \delta_{\tau_i} u) \tag{2.27}$$

$$\tau_{i+1} = \tau_i + r$$

$$z(t) = (z_1(t),\ldots,z_n(t))' \in D \subseteq \Re^n$$

where $P : C^0([0,r];\Omega) \times L^\infty([0,r];U) \to \Re^n$ is the operator defined by (2.23).

**Corollary 2.11:** *Consider system (1.1) and assume that hypotheses (H1), (H2) hold. Consider the unique solution $(x(t), y(t), z(t)) \in D \times \Omega \times D$ of (1.1), (2.27) with arbitrary initial condition $(x_0, y_0, z_0) \in D \times \Omega \times D$ corresponding to arbitrary input $u \in L^\infty_{loc}(\Re_+; U)$. Then the solution $(x(t), y(t), z(t)) \in D \times \Omega \times D$ of (1.1), (2.27) satisfies:*

$$z(t) = x(t), \text{ for all } t \geq r \tag{2.28}$$

**Example 2.12:** Consider the system

$$\dot{x}(t) = a(y(t), u(t))x(t)$$
$$\dot{y}(t) = f(y(t), u(t)) + c(y(t))x(t) \qquad (2.27)$$
$$x(t) \in D, \ y(t) \in \Omega, \ u(t) \in U$$

where $D = \Re$ or $D = (0, +\infty)$, $\Omega \subseteq \Re$ is an open set, $U \subseteq \Re$ is a closed non-empty set, $a : \Omega \times U \to \Re$, $f : \Omega \times U \to \Re$, $c : \Omega \to \Re$ are continuously differentiable mappings satisfying the following hypothesis:

**(H3)** $\dfrac{dc}{dy}(y) f(y, u) \neq 0$, for all $(y, u) \in \Omega \times U$ with $c(y) = 0$

We also assume that for every $(x_0, y_0) \in D \times \Omega$ and $u \in L^\infty_{loc}(\Re_+; U)$ there exists a unique solution $[0, +\infty) \ni t \to (x(t), y(t)) = (x(t, x_0, y_0; u), y(t, x_0, y_0; u)) \in D \times \Omega$ for (2.27) satisfying (2.27) for almost every $t \geq 0$. It is clear from (2.27), the fact that $D = \Re$ or $D = (0, +\infty)$ and the integral formula $x(t) = x(0) \exp\left(\int_0^t a(y(s), u(s)) ds\right)$ that hypothesis (H2) holds. Moreover, hypothesis (H3) and Corollary 2.4 guarantee that system (2.27) is strongly observable in time $r > 0$ (notice that $\det\left([C'(t, x_0, y_0; u) \Phi(t, x_0, y_0; u)]\right) = c(y(t, x_0, y_0; u)) \exp\left(\int_0^t a(y(s, x_0, y_0; u), u(s)) ds\right)$).

Using the formulas $\Phi(t, y, u) = \exp\left(\int_0^t a(y(s), u(s)) ds\right)$, $q(\tau) = \int_0^\tau c(y(s)) \exp\left(\int_0^s a(y(w), u(w)) dw\right) ds$,

$p(\tau) = y(\tau) - y(0) - \int_0^\tau f(y(s), u(s)) ds$, $Q = \int_0^r \left(\int_0^\tau c(y(s)) \exp\left(\int_0^s a(y(w), u(w)) dw\right) ds\right)^2 d\tau$, we conclude that the system:

$$\dot{z}(t) = a(y(t), u(t)) z(t), \ t \in [\tau_i, \tau_{i+1})$$
$$z(\tau_{i+1}) = P(\delta_{\tau_i} y, \delta_{\tau_i} u) \qquad (2.28)$$
$$\tau_{i+1} = \tau_i + r$$

where

$$P(y, u) = \exp\left(\int_0^r a(y(s), u(s)) ds\right) \frac{\int_0^r \left(y(\tau) - y(0) - \int_0^\tau f(y(s), u(s)) ds\right)\left(\int_0^\tau c(y(s)) \exp\left(\int_0^s a(y(w), u(w)) dw\right) ds\right) d\tau}{\int_0^r \left(\int_0^\tau c(y(s)) \exp\left(\int_0^s a(y(w), u(w)) dw\right) ds\right)^2 d\tau}$$

is a dead-beat reduced order observer for system (2.27).  ◁

# 3. Applications

This section is devoted to the study of two important applications: the operation of a batch chemical reactor and the robust estimation of the frequency of a sinusoidal signal.

**Application 1: Batch Reactor**

The system to be studied is similar to the system studied in [27]. We assume a batch reactor where the chemical reactions $A \to B \to C$ are taking place. The goal is the maximization of the concentration of the intermediate $B$. Both reactions are assumed to be exothermic and we also assume that the reactor is surrounded by jacket containing (heating or cooling) medium of constant temperature $T_s$ (see [19]). Finally, assuming that the kinetics of the reactions are first order, then the system is described by the following system of differential equations (see [19]):

$$\dot{c}_A = -k_1 \exp\left(-\frac{E_1}{T}\right) c_A$$
$$\dot{c}_B = k_1 \exp\left(-\frac{E_1}{T}\right) c_A - k_2 \exp\left(-\frac{E_2}{T}\right) c_B \quad (3.1)$$
$$\dot{T} = J_1 k_1 \exp\left(-\frac{E_1}{T}\right) c_A + J_2 k_2 \exp\left(-\frac{E_2}{T}\right) c_B + h(T_s - T)$$

where $c_A, c_B$ are the concentrations of $A, B$, respectively, $T$ is the temperature of the reactor, $J_1, J_2, k_1, k_2, h, E_1, E_2 > 0$ are positive parameters. The temperature is continuously measured. However, it is often very difficult to obtain accurate measurements of the concentrations of $A, B$. An observer for the estimation of $c_A, c_B$ is needed because (ideally) the operation of the reactor will be stopped when $c_B(t)$ is maximized, i.e., when $k_1 \exp\left(-\frac{E_1}{T(t)}\right) c_A(t) = k_2 \exp\left(-\frac{E_2}{T(t)}\right) c_B(t)$. Therefore, an accurate estimation of $c_A, c_B$ is needed. Furthermore, the observer must provide accurate estimations in short time. Indeed, if the system is left to operate for large times then the goal will not be achieved (because, we have $\lim_{t \to +\infty} c_A(t) = \lim_{t \to +\infty} c_B(t) = 0$).

System (3.1) is of the form (1.1) with $n = 2$, $y = T$, $x = (c_A, c_B)$, $A(y, u) = \begin{bmatrix} -k_1 \exp\left(-\frac{E_1}{y}\right) & 0 \\ k_1 \exp\left(-\frac{E_1}{y}\right) & -k_2 \exp\left(-\frac{E_2}{y}\right) \end{bmatrix}$,

$f(y, u) = h(T_s - y)$, $c_1(y) = J_1 k_1 \exp\left(-\frac{E_1}{y}\right)$, $c_2(y) = J_2 k_2 \exp\left(-\frac{E_2}{y}\right)$. Moreover, $O = D \times \Omega$, $D = (0, c_1) \times (0, c_2)$ and $\Omega = (T_{\min}, T_{\max})$, where $c_1 > 0$, $0 < T_{\min} \le T_s$, $\frac{k_1}{k_2} c_1 < c_2$ if $E_1 \ge E_2$ and $\frac{k_1}{k_2} \exp\left(\frac{E_2 - E_1}{T_{\min}}\right) c_1 < c_2$ if $E_1 < E_2$, $\frac{J_1 k_1 c_1 + J_2 k_2 c_2}{h} + T_s \le T_{\max}$ and $u$ and $U \subseteq \mathfrak{R}$ are irrelevant (because there is no input). The reader can easily verify that $D \times \Omega$ is positively invariant for (3.1) and that for every $(x_0, y_0) \in D \times \Omega$ there exists a unique solution $[0, +\infty) \ni t \to (x(t), y(t)) = (x(t, x_0, y_0; u), y(t, x_0, y_0; u)) \in D \times \Omega$ for (3.1) satisfying (3.1) for every $t \ge 0$. Moreover, hypothesis (H2) holds with

$$\Phi(t, y) = \begin{bmatrix} \exp\left(-\int_0^t k_1 \exp\left(-\frac{E_1}{y(s)}\right) ds\right) & 0 \\ k_1 \int_0^t \exp\left(-\frac{E_1}{y(\tau)} - \int_\tau^t k_2 \exp\left(\frac{-E_2}{y(s)}\right) ds - \int_0^\tau k_1 \exp\left(\frac{-E_1}{y(s)}\right) ds\right) d\tau & \exp\left(-\int_0^t k_2 \exp\left(-\frac{E_2}{y(s)}\right) ds\right) \end{bmatrix}$$

Next assume that one of the following hypotheses holds:

**(A1)** If $E_1 = E_2$ then $(J_1 + J_2)k_2 \neq J_1 k_1$. If $E_1 \neq E_2$ then there exists $a > 0$ such that for all $T \in (T_{\min}, T_{\max})$ with
$$\frac{T^2}{(E_2 - E_1)J_1 h} \exp\left(\frac{-E_2}{T}\right)\left[(J_1 + J_2)k_2 - J_1 k_1 \exp\left(\frac{E_2 - E_1}{T}\right)\right] + T > T_s \text{ it holds that:}$$

$$\frac{T^2}{(E_2 - E_1)J_1} \exp\left(\frac{-E_2}{T}\right)\left[(J_1 + J_2)k_2 - J_1 k_1 \exp\left(\frac{E_2 - E_1}{T}\right)\right] \leq -a \tag{3.2}$$

or

**(A2)** If $E_1 = E_2$ then $(J_1 + J_2)k_2 \neq J_1 k_1$. If $E_1 \neq E_2$ then there exists $a > 0$ such that for all $T \in (T_{\min}, T_{\max})$ with
$$\frac{T^2}{(E_2 - E_1)J_1 h} \exp\left(\frac{-E_2}{T}\right)\left[(J_1 + J_2)k_2 - J_1 k_1 \exp\left(\frac{E_2 - E_1}{T}\right)\right] + T > T_s \text{ it holds that:}$$

$$\frac{T^2}{(E_2 - E_1)J_1} \exp\left(\frac{-E_2}{T}\right)\left[(J_1 + J_2)k_2 - J_1 k_1 \exp\left(\frac{E_2 - E_1}{T}\right)\right] \geq a \tag{3.3}$$

We will show that system (3.1) under hypothesis (A1) or under hypothesis (A2) satisfies hypothesis (H1) as well for $r \geq \frac{T_{\max} - T_{\min}}{a}$ if $E_1 \neq E_2$ and for every $r > 0$ if $E_1 = E_2$. Therefore the system:

$$\begin{aligned}
\dot{z}_1(t) &= -k_1 \exp\left(-\frac{E_1}{T(t)}\right) z_1(t) \\
\dot{z}_2(t) &= k_1 \exp\left(-\frac{E_1}{T(t)}\right) z_1(t) - k_2 \exp\left(-\frac{E_2}{T(t)}\right) z_2(t) \quad , \quad t \in [\tau_i, \tau_{i+1}) \\
z(\tau_{i+1}) &= \Phi(r, \delta_{\tau_i} T) G(\delta_{\tau_i} T) \\
\tau_{i+1} &= \tau_i + r
\end{aligned} \tag{3.4}$$

where $G(T) = (G_1(T), G_2(T))$ for all $T \in C^0([0, r]; (T_{\min}, T_{\max}))$ is defined by

$$G_1(T) = \frac{\int_0^r \phi_2^2(t)dt \int_0^r \left(T(t) - T(0) - \int_0^t h(T_s - T(s))ds\right)\phi_1(t)dt - \int_0^r \phi_1(t)\phi_2(t)dt \int_0^r \left(T(t) - T(0) - \int_0^t h(T_s - T(s))ds\right)\phi_2(t)dt}{\left[\int_0^r \phi_1^2(t)dt \int_0^r \phi_2^2(t)dt - \left(\int_0^r \phi_1(t)\phi_2(t)dt\right)^2\right]}$$

$$G_2(T) = \frac{\int_0^r \phi_1^2(t)dt \int_0^r \left(T(t) - T(0) - \int_0^t h(T_s - T(s))ds\right)\phi_2(t)dt - \int_0^r \phi_1(t)\phi_2(t)dt \int_0^r \left(T(t) - T(0) - \int_0^t h(T_s - T(s))ds\right)\phi_1(t)dt}{\left[\int_0^r \phi_1^2(t)dt \int_0^r \phi_2^2(t)dt - \left(\int_0^r \phi_1(t)\phi_2(t)dt\right)^2\right]}$$

$$\phi_1(t) = (J_1 + J_2)\left[1 - \exp\left(-\int_0^t k_1 \exp\left(-\frac{E_1}{T(s)}\right)ds\right)\right] - J_2 k_1 \int_0^t \exp\left(-\frac{E_1}{T(\tau)} - \int_\tau^t k_2 \exp\left(-\frac{E_2}{T(s)}\right)ds - \int_0^\tau k_1 \exp\left(-\frac{E_1}{T(s)}\right)ds\right)d\tau$$

$$\phi_2(t) = J_2\left[1 - \exp\left(-\int_0^t k_2 \exp\left(-\frac{E_2}{T(s)}\right)ds\right)\right]$$

will be a reduced-order dead-beat observer which satisfies $z(t) = (c_A(t), c_B(t))$ for all $t \geq r$.

Thus, we are left to show that hypothesis (H1) holds, i.e., that every input $u \in L^\infty([0,r];U)$ strongly distinguishes every state $(x_0, y_0) \in D \times \Omega$ in time $r > 0$. We proceed by contradiction, i.e., we assume that there exists $u \in L^\infty([0,r];U)$, $(x_0, y_0) \in D \times \Omega$ such that $u \in L^\infty([0,r];U)$ does not strongly distinguish the state $(x_0, y_0) \in D \times \Omega$ in time $r > 0$. By virtue of Remark 2.4 and applying (2.14) with $n = 2$, $t_1 = 0$ in conjunction with the fact that

$$C'(t, x_0, y_0; u)\Phi(t, x_0, y_0; u) = [\mu_1 \quad \mu_2]$$

$$\mu_1 = k_1 J_1 \exp\left(-\frac{E_1}{T(t)}\right) \exp\left(-\int_0^t k_1 \exp\left(-\frac{E_1}{T(s)}\right) ds\right)$$

$$+ J_2 k_2 k_1 \exp\left(-\frac{E_2}{T(t)}\right) \int_0^t \exp\left(-\frac{E_1}{T(\tau)} - \int_\tau^t k_2 \exp\left(-\frac{E_2}{T(s)}\right) ds - \int_0^\tau k_1 \exp\left(-\frac{E_1}{T(s)}\right) ds\right) d\tau$$

$$\mu_2 = J_2 k_2 \exp\left(-\frac{E_2}{T(t)}\right) \exp\left(-\int_0^t k_2 \exp\left(-\frac{E_2}{T(s)}\right) ds\right)$$

where $T(t)$ is the component of the solution of (3.1) with initial condition $(x_0, y_0) \in D \times \Omega$, we obtain for all $t \in [0, r]$:

$$\begin{aligned} &J_1 \exp\left(-\frac{E_1}{T(t)}\right) \exp\left(-\int_0^t k_1 \exp\left(-\frac{E_1}{T(s)}\right) ds\right) \\ &+ J_2 k_2 \exp\left(-\frac{E_2}{T(t)}\right) \int_0^t \exp\left(-\frac{E_1}{T(\tau)} - \int_\tau^t k_2 \exp\left(-\frac{E_2}{T(s)}\right) ds - \int_0^\tau k_1 \exp\left(-\frac{E_1}{T(s)}\right) ds\right) d\tau = \\ &M \exp\left(-\frac{E_2}{T(t)}\right) \exp\left(-\int_0^t k_2 \exp\left(-\frac{E_2}{T(s)}\right) ds\right) \end{aligned} \quad (3.5)$$

where

$$M = J_1 \exp\left(\frac{E_2 - E_1}{y_0}\right) \quad (3.6)$$

By differentiating (3.5), we obtain for all $t \in [0, r]$:

$$(J_1 + J_2)k_2 \exp\left(-\frac{E_2}{T(t)}\right) - J_1 k_1 \exp\left(-\frac{E_1}{T(t)}\right) = (E_2 - E_1) J_1 \frac{\dot{T}(t)}{T^2(t)} \quad (3.7)$$

We next distinguish the cases:

<u>Case 1:</u> $E_1 = E_2$

In this case (3.7) implies $(J_1 + J_2)k_2 = J_1 k_1$, which leads to contradiction with hypothesis (A1) or hypothesis (A2). Thus, (3.5) cannot hold for all $t \in [0, r]$. Consequently, in this case we can conclude that hypothesis (H1) holds for every $r > 0$.

<u>Case 2:</u> $E_1 \neq E_2$

In this case (3.7) and (3.1) (which implies that $\dot{T} > h(T_s - T)$) imply that the following inequality must hold for all $t \in [0, r]$:

$$\frac{T^2(t)}{(E_2-E_1)J_1 h}\exp\left(\frac{-E_2}{T(t)}\right)\left[(J_1+J_2)k_2 - J_1 k_1 \exp\left(\frac{E_2-E_1}{T(t)}\right)\right] + T(t) > T_s \qquad (3.8)$$

We next assume that hypothesis (A1) holds. The analysis is completely similar for the case where hypothesis (A2) holds. Inequality (3.2) implies that there exists $a > 0$ such that

$$\dot{T}(t) \leq -a, \text{ for all } t \in [0,r] \qquad (3.9)$$

For every $T(0) = y_0 \in \Omega = (T_{min}, T_{max})$ the differential inequality (3.9) in conjunction with the fact that $r \geq \frac{T_{max} - T_{min}}{a}$ leads to the contradiction $T(r) \leq T_{min}$, i.e., $T(r) \notin \Omega$. Thus, (3.5) cannot hold for all $t \in [0,r]$. Consequently, in this case we can conclude that hypothesis (H1) holds for every $r \geq \frac{T_{max} - T_{min}}{a}$.

**Remark:** It should be emphasized that if $E_1 = E_2$ then the condition for strong observability of system (3.1) is $(J_1 + J_2)k_2 \neq J_1 k_1$. Indeed, notice that if $(J_1 + J_2)k_2 = J_1 k_1$ then system (3.1) gives the observable subsystem

$$\frac{d}{dt}\left[(J_1+J_2)c_A + J_2 c_B\right] = -\{(J_1+J_2)c_A + J_2 c_B\}k_2 \exp\left(-\frac{E}{RT}\right)$$
$$\dot{T} = k_2 \exp\left(-\frac{E}{RT}\right)\left[(J_1+J_2)c_A + J_2 c_B\right] + h(T_s - T) \qquad (3.10)$$

from which we can conclude that the measurement of the temperature can only give us estimates for the quantity $(J_1 + J_2)c_A + J_2 c_B$ and not for the components of the state vector $c_A, c_B$. Moreover, it should be noted that if $E_1 < E_2$ and $(J_1 + J_2)k_2 < J_1 k_1$, then there exists $a > 0$ such that hypothesis (A1) holds automatically. Similarly, if $E_1 > E_2$ and $(J_1 + J_2)k_2 > J_1 k_1$, then there exists $a > 0$ such that hypothesis (A2) holds automatically.

**Application 2: Estimation of Frequency of a Sinusoidal Signal**

The problem studied is to estimate the frequency $\omega > 0$ of a sinusoidal signal $y(t) = A\sin(\omega t + \varphi)$. The problem can be cast as an observer problem for the following system:

$$\begin{aligned}\dot{y}(t) &= x_1(t) \\ \dot{x}_1(t) &= x_2(t)y(t) \\ \dot{x}_2(t) &= 0 \\ (x_1(t), x_2(t), y(t)) &\in O \end{aligned} \qquad (3.11)$$

where $O := \{(x_1, x_2, y) \in \Re^3 : y^2 + x_1^2 > 0, x_2 < 0\}$ and $x_2(t) = -\omega^2$. It should be noticed that system (3.11) is forward complete and satisfies hypothesis (H1) for every $r > 0$. Indeed, only initial states on the manifold $y = x_1 = 0$ can give identical output responses. The results of the previous section can be applied in order to give the hybrid full-order observer:

$$\begin{aligned}\dot{w}(t) &= z_1(t) \\ \dot{z}_1(t) &= z_2(t)w(t) \\ \dot{z}_2(t) &= 0 \\ (z_1(t), z_2(t), w(t)) &\in O \end{aligned} \text{, for } t \in [\tau_i, \tau_{i+1}) \qquad (3.12)$$

$$z_1(\tau_{i+1}) = \frac{3\left(\int_{\tau_i}^{\tau_{i+1}}\phi^2(t)dt - \int_{\tau_i}^{\tau_{i+1}}y(\tau)d\tau \int_{\tau_i}^{\tau_{i+1}}(t-\tau_i)\phi(t)dt\right)\int_{\tau_i}^{\tau_{i+1}}(y(t)-y(\tau_i))(t-\tau_i)dt}{r^3\int_{\tau_i}^{\tau_{i+1}}\phi^2(t)dt - 3\left(\int_{\tau_i}^{\tau_{i+1}}(t-\tau_i)\phi(t)dt\right)^2}$$

$$+ \frac{\left(r^3\int_{\tau_i}^{\tau_{i+1}}y(\tau)d\tau - 3\int_{\tau_i}^{\tau_{i+1}}(t-\tau_i)\phi(t)dt\right)\int_{\tau_i}^{\tau_{i+1}}(y(t)-y(\tau_i))\phi(t)dt}{r^3\int_{\tau_i}^{\tau_{i+1}}\phi^2(t)dt - 3\left(\int_{\tau_i}^{\tau_{i+1}}(t-\tau_i)\phi(t)dt\right)^2}$$

(3.13)

$$z_2(\tau_{i+1}) = \frac{-3\int_{\tau_i}^{\tau_{i+1}}(t-\tau_i)\phi(t)dt \int_{\tau_i}^{\tau_{i+1}}(y(t)-y(\tau_i))(t-\tau_i)dt + r^3\int_{\tau_i}^{\tau_{i+1}}(y(t)-y(\tau_i))\phi(t)dt}{r^3\int_{\tau_i}^{\tau_{i+1}}\phi^2(t)dt - 3\left(\int_{\tau_i}^{\tau_{i+1}}(t-\tau_i)\phi(t)dt\right)^2}$$

(3.14)

$$w(\tau_{i+1}) = y(\tau_{i+1})$$
$$\tau_{i+1} = \tau_i + r$$

(3.15)

with $\phi(t) = \int_{\tau_i}^{t}\left(\int_{\tau_i}^{s}y(l)dl\right)ds$. It should be noted that hypothesis (H2) does not hold for system (3.11). The frequency $\omega$ is estimated by means of the formula $\hat{\omega} = \sqrt{-z_2(t)}$.

We assume next that the measured signal is corrupted by high frequency noise, i.e., we assume that

$$y(t) = A\sin(\omega t + \varphi) + a\sin(ft)$$

(3.16)

Exactly the same test of robustness as the one in [26] is performed: the parameters are chosen to be $A = 2$, $a = 0.2$, $\omega = 3$. Three cases are considered for the frequency of the noise: $f = 10$, $f = 100$ and $f = 1000$.

The effectiveness of formula (3.14) with $r = 1$, $\tau_i = 0$ is shown in Figures 1,2,3 as a function of the phase angle $\varphi$. It is shown that the greatest estimation error is $6.6\%$, $1.3\%$ and $0.083\%$ for the cases $f = 10$, $f = 100$ and $f = 1000$, respectively. The accuracy of the estimation is similar to the one obtained in [26], where the steady state estimation error was $10\%$, $1\%$ and $0.1\%$ for the cases $f = 10$, $f = 100$ and $f = 1000$, respectively. It should be noted that the estimated frequency for the hybrid observer is provided only after $r = 1s$, while in [26] at least $5s$ are needed in order to obtain an accurate estimate for the frequency.

However, if larger values for $r > 0$ are used then the accuracy of the estimation can be increased significantly. Figure 4 shows the estimated frequency from formula (3.14) with $\tau_i = 0$ as a function of $r$ for the case $f = 10$. The phase angle was selected to be $\varphi = 1.9$: this is the value of the phase angle that the largest error of the estimation occurs (see also Figure 1). For $r = 3$ the estimation error is $0.066\%$, i.e., it is 100 times less than the error obtained for $r = 1$.

Finally, it should be noted that the full-order observer (3.12), (3.13), (3.14), (3.15) can be used for system (3.11) even if the open set $O$ is defined to be $O := \{(x_1, x_2, y) \in \Re^3 : y^2 + x_1^2 > 0\}$. This is the case studied in [33].

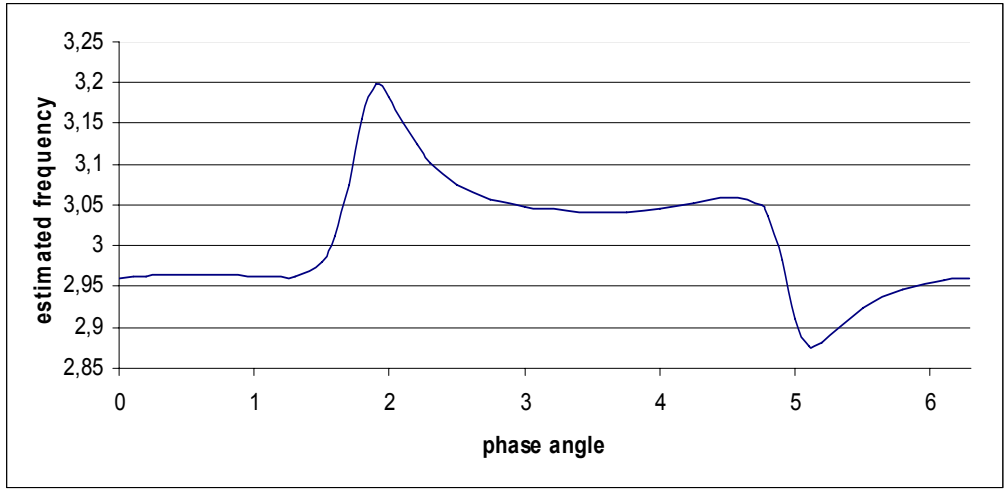

**Figure 1:** Estimated frequency from formula (3.14) with $r=1$, $\tau_i = 0$, $f=10$ as a function of the phase angle $\varphi$.

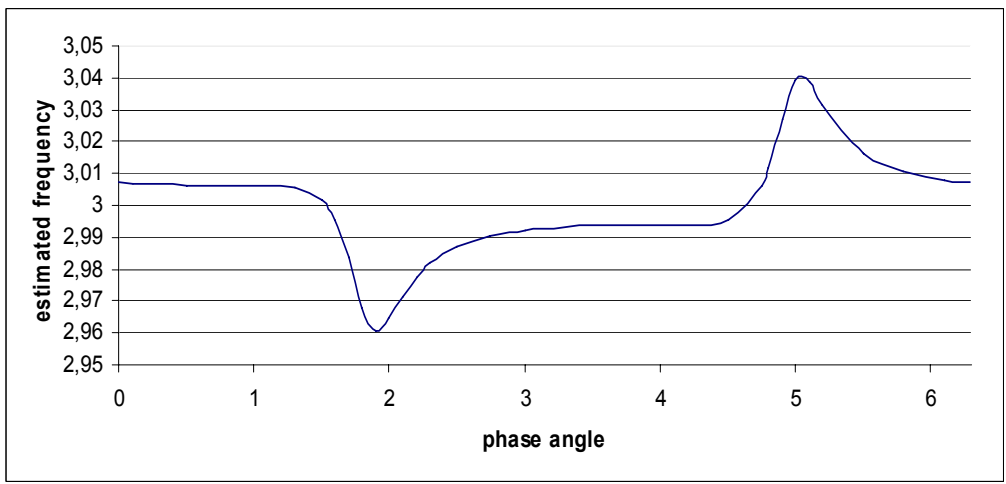

**Figure 2:** Estimated frequency from formula (3.14) with $r=1$, $\tau_i = 0$, $f=100$ as a function of the phase angle $\varphi$.

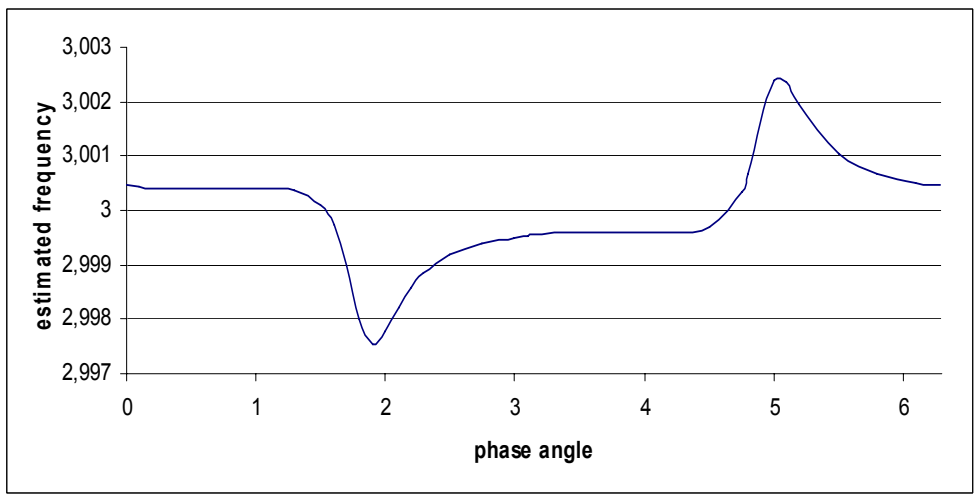

**Figure 3:** Estimated frequency from formula (3.14) with $r=1$, $\tau_i = 0$, $f=1000$ as a function of the phase angle $\varphi$.

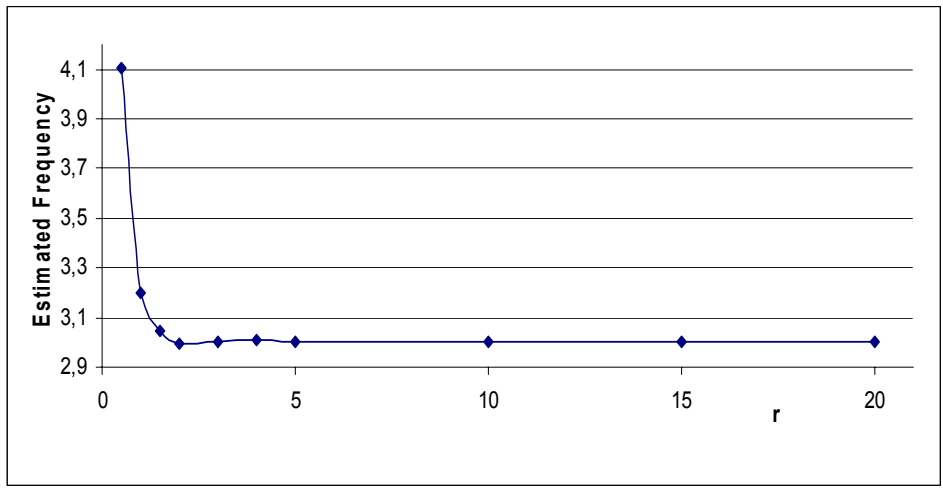

**Figure 4:** Estimated frequency from formula (3.14) with $\varphi = 1.9$, $\tau_i = 0$, $f = 10$ as a function of $r$.

## 4. Concluding Remarks

In this work, a novel hybrid strategy has been developed for solving the dead-beat observer design problem for a class of nonlinear systems with unmeasured states appearing linearly in the differential equations. To this end, the notion of strong observability of a nonlinear control system is introduced and utilized. The proposed methodology is applied to yield a hybrid full-order or reduced-order dead-beat observer for a batch reactor without concentration measurements. Moreover, the observer is used for the estimation of the frequency of a sinusoidal signal. The results show that accurate estimations can be provided even if the signal is corrupted by high frequency noise.

Future work can shed new light to the problem of dynamic output feedback stabilization, already studied in [3,8,9,22,23,24]. The dead-beat feature of the proposed observer implies that any static feedback stabilizer for (1.1) can be used in conjunction with the hybrid dead-beat observer (2.25), provided that the inputs produced by the applied feedback can distinguish all states in finite time and that the solution does not blow up during the initial transient period. Another direction for future work is the application of the hybrid, dead-beat observer to systems of mathematical biology: the chemostat model (see [29]) takes the form of system (1.1), when the nutrient concentration is measured.

**Acknowledgements:** This work has been supported in part by the NSF grants DMS-0504462 and DMS-0906659.

# Appendix

**Proof of Proposition 2.3:** We prove the implications (a) $\Leftrightarrow$ (b), (b) $\Leftrightarrow$ (c) and (c) $\Leftrightarrow$ (d)

(a) $\Rightarrow$ (b) The proof of this implication will be made by contradiction. Suppose that the input $u \in L^\infty([0,r];U)$ strongly distinguishes the state $(x_0, y_0) \in O$ in time $r > 0$. Notice that Fact I and definitions (2.6), (2.7), (2.8), (2.9) imply that problem (2.10) is always solvable and always admits the solution $\xi = x_0$ with

$$0 = \int_0^r |p(t, x_0, y_0; u) - q'(t, x_0, y_0; u) x_0|^2 dt = \min_{\xi \in B(y_0)} \int_0^r |p(t, x_0, y_0; u) - q'(t, x_0, y_0; u) \xi|^2 dt \quad \text{(A1)}$$

Consequently, the negation of (b) is the following statement:

"Problem (2.10) admits the solution $\xi = x_1 \in \Re^n$ with $x_1 \neq x_0$ and $(x_1, y_0) \in O$"

Therefore, we assume that the above statement holds. By virtue of (A1) we must have $0 = \int_0^r |p(t, x_0, y_0; u) - q'(t, x_0, y_0; u) x_1|^2 dt$. Continuity of the mappings $t \to p(t, x_0, y_0; u)$ and $t \to q(t, x_0, y_0; u)$ implies that the following statement holds:

"there exists $x_1 \neq x_0$ with $(x_1, y_0) \in O$ such that $p(t, x_0, y_0; u) = q'(t, x_0, y_0; u) x_1$ for all $t \in [0, r]$."

The above statement in conjunction with definitions (2.6), (2.7), (2.9) shows that (by direct differentiation):

$$\frac{d}{dt} y(t, x_0, y_0; u) = f(y(t, x_0, y_0; u), u(t)) + C'(t, x_0, y_0; u)(\Phi(t, x_0, y_0; u) x_1 + \theta(t, x_0, y_0; u))$$

where $f(y, u) := (f_1(y, u), \ldots, f_k(y, u))'$ and

$$\frac{d}{dt}\left((\Phi(t, x_0, y_0; u) x_1 + \theta(t, x_0, y_0; u))\right)$$
$$= A(y(t, x_0, y_0; u), u(t))\left((\Phi(t, x_0, y_0; u) x_1 + \theta(t, x_0, y_0; u))\right) + b(y(t, x_0, y_0; u), u(t))$$

for almost all $t \in [0, r]$. Consequently, uniqueness of solutions for (1.1) implies that $x(t, x_1, y_0; u) = \Phi(t, x_0, y_0; u) x_1$ and $y(t, x_1, y_0; u) = y(t, x_0, y_0; u)$ for all $t \in [0, r]$. Hence, it holds that:

"there exists $x_1 \neq x_0$ with $(x_1, y_0) \in O$ such that $y(t, x_1, y_0; u) = y(t, x_0, y_0; u)$ for all $t \in [0, r]$."

The above statement contradicts the assumption that the input $u \in L^\infty([0, r]; U)$ strongly distinguishes the state $(x_0, y_0) \in O$ in time $r > 0$.

(b) $\Rightarrow$ (a) Again the proof of this implication will be made by contradiction. Suppose that problem (2.10) admits the unique solution $\xi = x_0$.

Assume that the input $u \in L^\infty([0, r]; U)$ does not strongly distinguish the state $(x_0, y_0) \in O$ in time $r > 0$. This implies that

"there exists $x_1 \in \Re^n$ with $x_1 \neq x_0$ and $(x_1, y_0) \in O$ such that $y(t, x_1, y_0; u) = y(t, x_0, y_0; u)$ for all $t \in [0, r]$."

The reader should notice that the $y-$components of the different initial states $(x_0, y_0) \in O$ and $(x_1, y_0) \in O$ which produce identical outputs for $t \in [0, r]$, necessarily coincide. By virtue of Fact I and definitions (2.6), (2.7), (2.8), (2.9) it follows that $p(t, x_0, y_0; u) = q'(t, x_0, y_0; u) x_1$ for all $t \in [0, r]$. Hence, we must have:

$$0 = \int_0^r |p(t, x_0, y_0; u) - q'(t, x_0, y_0; u) x_1|^2 dt$$

The above equality shows that $x_1 \in \Re^n$ with $x_1 \neq x_0$ and $(x_1, y_0) \in O$ is a solution of problem (2.10), which contradicts the uniqueness of the solution for problem (2.10).

(b) $\Rightarrow$ (c) Again the proof of this implication will be made by contradiction. Suppose that problem (2.10) admits the unique solution $\xi = x_0$. Notice that the objective function for problem (2.10) is the quadratic function:

$$R(\xi) := \int_0^r \left| p(t, x_0, y_0; u) - q'(t, x_0, y_0; u) \xi \right|^2 dt$$

$$= \int_0^r \left| p(t, x_0, y_0; u) \right|^2 dt - 2 \int_0^r p'(t, x_0, y_0; u) q'(t, x_0, y_0; u) dt\, \xi + \xi' Q(r, x_0, y_0; u) \xi \tag{A2}$$

where $Q(r, x_0, y_0; u)$ is defined by (2.11) and for which it holds that

$$R(\xi) = R(x_0) + 2 \left[ -\int_0^r q(t, x_0, y_0; u) p(t, x_0, y_0; u) dt + Q(r, x_0, y_0; u) x_0 \right]' (\xi - x_0) + (\xi - x_0)' Q(r, x_0, y_0; u)(\xi - x_0) \tag{A3}$$

Since $\xi = x_0$ is a solution of problem (2.10) and since $O$ is open the above equality shows that we must necessarily have

$$\int_0^r q(t, x_0, y_0; u) p(t, x_0, y_0; u) dt = Q(r, x_0, y_0; u) x_0 \tag{A4}$$

On the other hand assume that statement (c) does not hold, i.e., assume that the symmetric and positive semidefinite matrix $Q(r, x_0, y_0; u) := \int_0^r q(t, x_0, y_0; u) q'(t, x_0, y_0; u) dt$ is not positive definite. Therefore there exists $\zeta \in \Re^n$, $\zeta \neq 0$ such that $0 = \zeta' Q(r, x_0, y_0; u) \zeta$. It follows from (A3), (A4) and for sufficiently small $\lambda > 0$ that the vector $\xi = x_0 + \lambda \zeta$ will satisfy $(\xi, y_0) \in O$ (because $O$ is open) and $R(\xi) = R(x_0)$, i.e., the vector $\xi = x_0 + \lambda \zeta$ is an additional solution of problem (2.10) with $\xi \neq x_0$, a contradiction.

Therefore $Q(r, x_0, y_0; u) := \int_0^r q(t, x_0, y_0; u) q'(t, x_0, y_0; u) dt$ is positive definite. Equation (2.12) is a direct consequence of equation (A4).

(c) $\Rightarrow$ (b) This implication is a direct consequence of (A2), (A3), (A4), which show that

$$R(\xi) = (\xi - x_0)' Q(r, x_0, y_0; u)(\xi - x_0), \text{ for all } \xi \in \Re^n \text{ with } (\xi, y_0) \in O$$

Notice that equality (2.5) guarantees that $R(x_0) = 0$.

(c) $\Rightarrow$ (d) This implication follows from the fact that $\xi' Q(r, x_0, y_0; u) \xi := \int_0^r \left| q'(t, x_0, y_0; u) \xi \right|^2 dt$, for all $\xi \in \Re^n$.

(d) $\Rightarrow$ (c) Statement (c) follows from the fact that $\xi' Q(r, x_0, y_0; u) \xi := \int_0^r \left| q'(t, x_0, y_0; u) \xi \right|^2 dt$, for all $\xi \in \Re^n$ and the fact that the mapping $t \to q(t, x_0, y_0; u)$ is continuous. Equality (2.12) is a direct consequence of Fact I and equality (2.5).

The proof is complete. ◁